\documentclass[runningheads]{llncs}
\usepackage[T1]{fontenc}
\usepackage{graphicx}
\usepackage{hyperref}
\usepackage{color}

\usepackage{amsmath}

\DeclareMathOperator{\expm1}{expm1}
\usepackage{amssymb}
\begin{document}
\title{Testing interval arithmetic libraries, including their IEEE-1788 compliance}
\titlerunning{Testing interval arithmetic libraries}
%
\author{Nathalie Revol\inst{1}\orcidID{0000-0002-2503-2274} 
\and Luis Benet\inst{2}\orcidID{0000-0002-8470-9054}
\and Luca Ferranti\inst{3}\orcidID{0000-0001-5588-0920}
\and Sergei Zhilin\inst{4}\orcidID{0000-0002-6961-1358}}
\authorrunning{N. Revol et al.}
%
\institute{INRIA - LIP UMR 5668, ENS Lyon, University Lyon 1, Inria, CNRS - Lyon, France
\email{Nathalie.Revol@inria.fr}\\
\and
Instituto de Ciencias F\'isicas, Universidad Nacional Aut\'onoma de M\'exico, Cuernavaca, Mexico
\email{benet@icf.unam.mx}\\
\and
University of Vaasa, Vaasa, Finland\\
\email{luca.ferranti@uwasa.fi}\\
\and
CSort LLC, Barnaul, Russia\\
\email{szhilin@gmail.com}}
\maketitle              
\begin{abstract}

As developers of libraries implementing interval arithmetic, we faced the same difficulties when it came to testing our libraries.
What must be tested? How can we devise relevant test cases for unit testing? How can we ensure a high (and possibly 100\%) test coverage?
In this paper we list the different aspects that, in our opinion, must be tested, giving indications on the choice of test cases.
Then we examine how several interval arithmetic libraries actually perform tests.
Next, we introduce two frameworks developed specifically to gather test cases and to incorporate easily new libraries in order to test them, namely JInterval and ITF1788.
Not every important aspects of our libraries fit in these frameworks and we list extra tests that we deem important, but not easy, to perform.

\end{abstract}

\section{Introduction}
\label{section-introduction}
Many libraries implement interval arithmetic, from XSC~\cite{XSC} and FiLib~\cite{FilibPub} among the pioneers, to Octave/Interval~\cite{OctaveIntervalPub}, ValidatedNumerics.jl~\cite{ValidatedNumericsCodes} and Moore~\cite{Mascarenhas2016} for the more recent ones.
In 2010, a comparison of several libraries~\cite{ICMS10} indicated that the underlying approach of the definition of interval arithmetic was, more or less, different and specific to each library.
It was thus impossible to compare the results obtained by different libraries, as so many differences impacted the computations.
To enable comparisons, it was decided to standardize interval arithmetic.
A collective effort, launched in 2008, led to the standardization of interval arithmetic, specified in the IEEE 1788-2015 standard~\cite{IEEE-1788-2015}.

The development of IEEE 1788-2015 compliant libraries, and thus the comparison of their results, was therefore made possible.
However, before even considering comparisons between libraries, the developers of an interval arithmetic library need to check their developments, to ensure that their library behaves correctly, in particular with respect to IEEE 1788-2015 compliance. How do they usually proceed? Through intensive tests. Formal proof that each operation is correct is the next step, but it will not be considered in what follows.

Another use of tests is, for users, at installation time, to check that the newly installed library behaves properly on one's platform. This is not a 100\% guarantee of correctness, but it increases the user's confidence.

Developers test for many aspects of their libraries, not only for IEEE 1788-2015 compliance, but also for instance the handling of exceptional or invalid inputs and so on. {\em Unit tests} are very common: for each operation or function, one prepares a set of inputs and checks whether the library computes the expected output, or at least an enclosure of it.
Rapidly, the developers realize that, in order to get sufficient confidence in the implementation of a given operation, in order to test a high enough coverage of the code (or 100\% of it), a huge number of cases must be tested. This implies that, {\em for each implemented operation}, a large part of the development time is devoted to devising these test cases.
We are developers of different libraries, we faced the same situation and we elaborated similar tests.
We decided to put our expertise in common and to share it. The
preliminary result of this collaboration is the content of this article.

In what follows, we propose to list, in Section~\ref{section-tests}, the different features that must be tested.
We then compare, in Section~\ref{section-libraries}, this list of recommended tests with the tests that are actually included in some well-known libraries (for which we have access to their test suites).
Our goal is to devise test cases for each situation we have identified and to share them:
in Section~\ref{section-frameworks} we introduce two existing frameworks that, on the one hand, ease the integration of a new library to be tested,
and that, on the other hand, offer a significant set of test cases, as well as a convenient mechanism to add new ones.
As no setup is perfect, in Section~\ref{section-wishlist}, we then give a list of other types of tests deemed useful in our development of interval arithmetic libraries, but still missing.
Our short- to midterm-goal is to enrich 
and share
test suites, and enhancing the capabilities of the frameworks introduced in Section~\ref{section-frameworks}.

\section{What must be tested?}
\label{section-tests}

Testing floating-point arithmetic has been studied thoroughly, as testified by Nelson Beebe's web page for software~\cite{Beebe} 
and more generally by the huge amount of bibliographic references and of programs available. 
Similarly, in this article, we gather and organize the important features that must be  tested in order to assess the correctness, compliance and possibly quality of a library implementing interval arithmetic.

\subsection{General remarks about unit tests}
Let us denote by \(\cal{L}\) the tested library, and intervals using boldface: \(\mathbf{x}\), \(\mathbf{y}\).
We focus here on unit tests, for a function, denoted by \(f\).
A unit test case is a pair composed of the input argument \(\mathbf{x}\) and the expected output \(\mathbf{y}\).

First, the output \(\mathbf{y}\) must be the tightest representable interval enclosing \(f (\mathbf{x}) \); otherwise a very accurate library could compute \(\mathbf{z}\) such that \( \mathbf{z} \subsetneq \mathbf{y} \) and still \( \mathbf{z} \supset f(\mathbf{x}) \).
To ensure the tightness property for \(\mathbf{y}\),
typically, one computes the endpoints \(\underline{y}\) and \(\bar{y}\) of \(\mathbf{y}=[\underline{y}, \bar{y}]\) using a precision higher than the computing precision of \(\cal{L}\). Let us illustrate this with \(\cal{L}\) using {\tt Binary64} floating-point numbers, on the (simple) example of the exponential function. We assume that the given floating-point implementation of \(\exp\) does not provide correct rounding, but that, for any precision \(q\), it satisfies \(RD_q(\exp(x)) \leq \exp(x) \leq RU_q (\exp (x)) \) for \(RD_q\) rounding downwards and \(RU_q\) rounding upwards, both in precision \(q\).
Given \(\mathbf{x} = [ \underline{x}, \bar{x}]\), to compute the infimum \(\underline{y}\) of \(\mathbf{y}\), one gets the approximations \(RD_q(\exp(\underline{x}))\) and \( RU_q(\exp(\underline{x}))\) of \(\exp(\underline{x})\) in high precision \(q\), and finally round them downwards in the target precision \(=53\). If \(RD_{53}(RD_q(\exp(\underline{x})))\) and \(RD_{53}(RU_q(\exp(\underline{x})))\) are equal, then they are the sought value \(\underline{y}\) for the infimum of \(\exp(\mathbf{x})\).

Once the test cases are devised, what conclusions can be drawn from the comparison between \(\mathbf{y}\) and \(\mathbf{z}\) the result computed by \(\cal{L}\)? Requiring equality may be too demanding. Inclusion is required, but we want to dismiss, for instance, an implementation of \(\sin\) that returns \([-1,1]\) for any argument. If the accuracy of the function is given in the specification of the library, how can it be used and checked?

A last general recommendation is to incorporate the following procedure: if \( (\mathbf{x}, \mathbf{y})\) is given as a pair of input and output, pick at random (a reasonably large number of) values \(x \in \mathbf{x}\) and check whether \(f(x) \in \mathbf{y}\) with \(f(x)\) computed by the underlying arithmetic, and whether \(f([x,x]) \subset \mathbf{y}\). Those tests are disconnected from the knowledge of the implementation of \(f\) and may hit a zone not considered (forgotten) in the development of \(f\), such as an overlooked quadrant for a trigonometric function.

\subsection{Tests common to all interval arithmetic libraries}

\begin{itemize}
\item {\bf Easy test cases.}
First, easy cases are tested: these test cases are both easy to devise and easy to compute. They constitute preliminary tests, at an early stage of the implementation of the function, to identify and correct the most obvious bugs. These test cases are chosen to cover, roughly speaking, the various possible magnitudes of the arguments, but without any exhaustivity.
They also span positive and negative values, if the domain and range of the function so permit.

\item {\bf Special and exceptional values.}
A time-consuming task, when implementing a function, is the handling of special and exceptional values.
In floating-point arithmetic, \(0\) is a special value, because it has two representations, \(+0\) and \(-0\), that is, signed zeros.
It is thus valuable to test whether one gets the same result, independently of the sign of the representation of \(0\), when \(0\) is an endpoint of either an input argument or of the output.

\(1\) is not a special value for floating-point arithmetic. However, it is often considered as a special value, for instance for the logarithmic functions and the arc-cosine: its image is the special value \(0\).

Infinities are also special values in floating-point arithmetic. If they are supported by the library, one must have test cases containing infinities as their endpoints, both for inputs and for outputs. If infinities are not supported by the library, it is even the more so important to test input arguments that yield overflow, such as the addition or the multiplication of
\([0, \mbox{\small MAXREAL}]\) by \([0, \mbox{\small MAXREAL}]\) or \([0.75*\mbox{\small MAXREAL}, \mbox{\small MAXREAL}]\) by \([0.75*\mbox{\small MAXREAL}, \mbox{\small MAXREAL}]\).
Another example is \( \exp[\underline{x}, \bar{x}])\) with \(\underline{x}\) such that \(\exp(\underline{x})\) rather small and \(\bar{x}\) such that \(\exp(\bar{x})\) overflows, or both \(\underline{x}\) and \(\bar{x}\) such that \(\exp(\bar{x})\) overflow.

Again, when floating-point arithmetic is used, {\tt NaN} is an exceptional value, that is the result of an invalid computation. 
How is a {\tt NaN} handled when it occurs as the endpoint of an input interval? The answer should be that, even if it has no mathematical meaning, in practice {\tt NaN}s propagate.
Some libraries use {\tt NaN} to denote the emptyset, as is the case in Intlab~\cite{Rump1999}.
Care must then be taken for the union or intersection of two intervals.

\item {\bf Cornercases.}
Another family of test cases is designed in order to test the difficulties encountered during the implementation, such as difficult-to-round values, or values at the boundaries of the domain of the function, or close to points whose image is an extremum.
A typical example is a floating-point input close to an integer multiple of \(\frac{\pi}{2}\), for trigonometric functions.
The various tools for testing {\em floating-point} arithmetic can be a source of inspiration. 
Indeed, with floating-point arithmetic, as there is a sudden change in absolute error between two consecutive floating-point numbers when the binade changes, it is worth exploring several binades, both for the inputs and for the output values. This is particularly true for mid-rad representation of the intervals.
With floating-point arithmetic, another delicate zone concerns subnormal numbers. Test cases should contain intervals 
\([0,\bar{x}]\) and \([\underline{x}, \bar{x}]\) such that \(\underline{x}\) and \(\bar{x}\) are subnormals -- for the logarithmic functions, or intervals \([\underline{x}, \bar{x}]\) such that their exponential has one or two subnormal endpoints.

\item {\bf Functions specific to intervals.}
Functions such as the union (convex hull of the union) or the intersection must be tested when one, or both endpoints  are special or exceptional values, or when one argument is the empty set, if it is supported by the library.

Functions such as the midpoint or radius have been thoroughly studied in~\cite{Goualard2014}, test cases can be taken from this reference.
\item {\bf Input and Output.}
Finally, I/O is the place where the most unexpected things can happen: much creativity is needed to cover a large variety of input values given as strings.
\end{itemize}

\subsection{Tests about IEEE-1788 compliance}
The IEEE 1788-2015 standard~\cite{IEEE-1788-2015} was developed to enable comparisons of interval methods and their implementation.
It mandates operations in Section 9;
it provides "hooks" to integrate several flavors of interval arithmetic and it defines the so-called "set-based" flavor along with its set of recommended operations in Section 10;
it mandates that intervals are decorated with flags that sum up the history of their computation (such as "every operation and function involved in this computation is defined and continuous on its arguments") in Section 8. It also specifies the handling of exceptions in Section 12.1.3 and the grammar for I/O in Section 13. Last but not least, it offers the possibility to return either the tightest possible result, or simply an accurate one, or a valid one in Section 12.10.1. "Compressed arithmetic" will not be mentioned, as it is not yet widely used.
Testing IEEE 1788-2015 compliance means testing all of these aspects:
\begin{itemize}
\item {\bf testing flavor compliance:} each
test case must be accompanied with the indication of the corresponding flavor: how can it be specified? If several flavors are tested (once they are defined by a new revision of the standard), 
either the testing framework or the library must offer a mechanism to change the current flavor in use.This is needed to test for instance \(\emptyset + [1,2]\) or \(\exp([-\infty,0])\) in the set-based flavor and then \([3,2]+[4,1]\) in the (future?) Kaucher flavor.
\item {\bf required functions / recommended functions:} the
standard mandates that some functions are implemented by the library, such as \(\cosh\), but only recommends other ones, such as \(\expm1\) for the set-based flavor. The tests must allow for executing functions if they exist in the library but not crashing when a recommended function is not available.
\item {\bf testing decorations:} the
set-based flavor defines 5 decorations. To test every possible meaningful combination of decorations for inputs and output, one must devise close to 25 test cases for unary operations and up to 125 test cases for binary operations. Care must be taken to exhaust the complete list of possibilities.
\item {\bf testing the different accuracies: tightness/accurate/valid:} we
already discussed the difficulty to validate a computed result, to compare it with the expected output, when it is not the tightest possible one.
The standard defines three possible levels of accuracy.
Test cases must include all three possibilities for each pair of input and output, and check if the corresponding levels are available in the library.
They must include the description of this level, the testing framework must be able to check whether the result is tight or valid (both are easy), or accurate (by computing the enclosing interval that is only slightly larger, as defined by the standard).
\item {\bf exceptional behavior:} exceptions
have already been mentioned above. The standard defines the different exceptions that must be signalled, the test cases must include these exceptions and check whether the correct exception is signalled.
\item {\bf standard specific I/O system:} the
standard defines precisely the various forms of inputs and outputs. At least one test case per form must be present in the test case. If one wants to test every possibility, with at most 3 values for each variable field lengths, one obtains a thousand test cases.
\end{itemize}

\subsection{Tests specific to some libraries}
Some libraries must need specific test cases, corresponding to their specificities.

\begin{itemize}
\item {\bf Representations of intervals.}
For instance, a library that uses the mid-rad representation of intervals must test that operations and I/O operate correctly on this representation.
There are libraries that keep track of the openness or closedness of each endpoint: test cases must be provided to check that this information is correctly computed. If each possibility is tested, this can multiply by 4 the number of test cases.
Some libraries do not employ floating-point types to represent intervals, either in inf-sup or mid-rad representation. Exact rational numbers can be used instead, as in MPRIA \cite{MPRIA} and JInterval~\cite{JIntervalPub}.
Other libraries also allow complex intervals, that is, they use complex arithmetic as their underlying arithmetic. 
\item {\bf Precision.}
A library that accomodates several floating-point formats, such as {\tt Binary32} and {\tt Binary64}, must test that it operates soundly with both formats and with combination of both. Moreover, a library that uses arbitrary precision, such as MPFI~\cite{MPFI} or ARB~\cite{ARB} must include tests for largely varying precisions, without claiming exhaustivity.
\item {\bf Directed roundings.}
In order to guarantee the inclusion property, libraries rely on directed roundings. They must have tests which are specific to how directed rounding is achieved.
		If it is achieved by changing the rounding mode, ideally it should be tested that this is thread safe and that using multiple threads doesn't lead to undefined behavior.
		If it is achieved via software, by emulating directed rounding in round to nearest (as it is done e.g. by default in IntervalArithmetic.jl), then they very likely rely heavily on the use of Error-Free Transforms, or EFT in short. Thus, hard to round cases and corner cases for EFT (when it overflows, underflows, returns NaN) should also be tested.
\end{itemize}

\subsection{On tests timing}
An important issue related to tests is time. To illustrate this question, let us mention that MPFI currently implements about 30 test cases for each function, and running the tests ({\tt make check}) takes 20 seconds of user's time on a reasonably fast PC.
If all these recommended test cases were implemented, there would be at least a few hundred test cases, at most a few thousands. The testing time will be multiplied according by a factor between 10 and 100. Shall we stick to our limited number of test cases and risk to have not enough tests performed, or shall we risk to have a more complete test cases coverage and no test performed at all, because of the time it takes?

\subsection{Need for a unified framework for testing interval arithmetic libraries}

While there are several shared features between all interval arithmetic libraries, they are all also somehow unique in design choices (how to handle direct rounding, tightness vs speed, etc.). This call for a unified framework to test interval arithmetic libraries. This has several benefits. First, a standard framework makes it easier to compare libraries against each other. It would also offer a tool for developers, to easily verify whether their library is compliant with the standard or not. Finally, an important use case of interval arithmetic is to perform rigorous computations, that can be used as mathematical proofs. While several theorems for rigorous computations have been proved, they all rely on the assumption that the underlying implementation of interval arithmetic is correct. For this reason, a unified testing framework, developed through inputs from several interval libraries developers, would increase the reliability of \textit{all} libraries.

Based on the above discussion, we identify the following requirements in a unified interval testing framework:
        \begin{enumerate}
                \item \textbf{Modularity:} tests should be structured in a modular way, to allow users and developers choosing what parts to test. For example, tests for decorated intervals should be separated by tests for bare intervals. Test for recommended functions should be separated by tests for required functions.
                \item \textbf{Completeness:} there should be tests for all functions required and recommended by the standards. The tests should contain both some simple normal cases and cornercases, such as exceptional behavior (handle NaN, overflow, underflow) and hard to round cases for finite precision data types.
                \item \textbf{Support for different number systems}: allow to test when the bounds of the interval are {\tt Binary64}, {\tt Binary32}, arbitrary precision floats etc.
                \item \textbf{Test for tightness:} verify that the computed interval is the smallest valid interval.
                \item \textbf{Test for validity:} achieving tight bounds can be challenging or at least computationally expensive. Libraries at early stages or libraries with a focus on efficiency may prefer to test only that the computed interval is valid, i.e. an enclosure of the expected result. A unified framework should allow the user to switch between testing tightness and validity, ideally also allowing to set a threshold on how much wider the returned interval is allowed to be.
                \item \textbf{Self-validated}: it should be possible for the user to verify that the results expected by the tests are indeed correct.
        \end{enumerate}

\section{Interval arithmetic libraries and their test sets} 
\label{section-libraries}
During decades, several interval arithmetic libraries have been developed, each with its own design choices. Since we are interesting in testing, here we only list some of the ones for which tests are freely available.

\begin{enumerate} \setlength\itemsep{1em}

    \item \textbf{MPFI} \cite{MPFI}: C library for arbitrary precision interval arithmetic. Based on MPFR.
Unit tests for each function: some "easy" values, exceptions, exact cases. Currently not IEEE-1788 compliant.

	\item \textbf{JInterval} \cite{JIntervalPub,JIntervalCodes}: Java library providing preliminary standard IEEE P1788 compliant exact implementation of interval arithmetic operations and controllable arbitrary precision for elementary functions. Use the set of rational numbers extended by $-\infty$ and $+\infty$ as a basic number type for intervals endpoints. On top of this implementation, tightest approximations of interval arithmetic and elementary functions are supported for float types \texttt{binary16}, \texttt{binary32}, \texttt{binary64}, \texttt{binary128}, etc. 

    \item \textbf{libieeep1788} \cite{libieeep1788Pub,libieeep1788Codes}:  C++ implementation of the preliminary IEEE P1788 standard. Relies on MPFR. Unit tests for IEEE-1788 compliance.

    \item \textbf{Octave interval package} \cite{OctaveIntervalPub,OctaveIntervalCodes}: GNU Octave interval arithmetic package compliant with the standard. Test sets derived from \texttt{libieeep1788, MPFI, C-XSC, FILIB} plus several new tests for most of the standard functions. These tests constitute the current ITF1788 testsuite. Relies on MPFR for both arithmetic operations and elementary functions.

    \item \textbf{ValidatedNumerics.jl} \cite{ValidatedNumericsPub,ValidatedNumericsCodes}: Julia set of packages for interval arithmetic and applications. Currently not conformant with the standard. By default, direct rounding is handled via the software emulator \textsf{RoundingEmulator.jl}. However, it also supports changing rounding modes. Correctly rounded elementary functions are computed via CRlibm when possible, and using MPFR as fallback. For {\tt Binary64} bounds, uses the ITF1788 testsuite. Also has hardcoded tests for other non-standard functionalities (complex intervals, interval boxes, other number formats)

	\item \textbf{GAOL} \cite{Goualard2005}: C++ library not compliant with the standard, e.g. lacks decorated intervals. Tight bounds on arithmetic operations using round upwards mode only. Supports CRlibm and  IBM Accurate Portable Mathematical Library for elementary functions computations. Uses its own unit tests.
	
\end{enumerate}

\section{Existing frameworks}
\label{section-frameworks}

\subsection{Why they are needed?}
\label{subsection-why-frameworks} 

As a motivating example, consider Kuliamin's paper\cite{Kuliamin2007} about testing floating-point mathematical functions. The paper defines several criteria for devising tests, such as considering exceptional cases (overflow, underflow, NaN), inputs out of domain, special values (signed zero, subnormal numbers) etc.

However, this leads to more than 15,000 values to test for the exponential function alone. When interval arithmetic is tested, this means that each endpoint must take 15,000 values: the combinatorics (even if one restricts the tests to proper intervals, where the left endpoint is not greater than the right endpoint) are in defavor of developing individually, manually, such intensive tests for each new library. 

A natural way to increase effectiveness of individual libraries developers is to unite their efforts in order to accumulate tests in a properly structured common shared database of tests equipped by  appropriate tools for tests usage. That is why the frameworks briefly described in the next two subsections are worth our attention.

\subsection{JInterval P1788 Test Launcher}

\texttt{P1788 Test Launcher} \cite{JIntervalPub,JIntervalLauncher} is an interval arithmetic libraries testing application based on the JInterval library. The features of JInterval entail a freedom for \texttt{P1788 Test Launcher} in the selection of computation modes corresponding to a tested library.

\texttt{P1788 Test Launcher} loads dynamic libraries (.so/.dll) with third party implementation of P1788 and checks the results obtained from a library with the tightest results computed internally using JInterval. Operations and functions from a tested library are called through the unified adaptor wrapper interface which must be preliminary implemented for the library being checked. The structure of the wrapper interface is adopted from C++ templates of \texttt{libieeep1788} library. So, technically, the implementation of the interface in most of cases boils down to overloading methods of \texttt{libieeep1788} classes describing set-based flavor of interval arithmetic for \texttt{Binary64} value set as a rule. Wrappers for the following interval arithmetic libraries are already included in \texttt{P1788 Test Launcher}: \texttt{Profil/BIAS}, \texttt{boost/interval}, \texttt{C-XSC}, \texttt{Filib}, \texttt{libieeep1788}, \texttt{libMoore}, \texttt{MPFI}. 

The launcher reads tests from plain text files of simple human-readable format and writes the results computed using tested library and JInterval and their relation into a plain text report. Test set included to the Launcher consists of over 14,000 tests which partly originated from \texttt{libieeep1788}, \texttt{Filib}, \texttt{libMoore} while the others are original.

\subsection{ITF1788 -- Interval Test Framework for IEEE 1788}

The idea behind ITF1788 is to simplify the development of unit tests for an interval arithmetic library independently of the programming language it is written in, the testing framework used and library peculiarities, such as custom function names. ITF1788 is a Python engine for converting pre-calculated tests from domain specific language, called the Interval Test Language (ITL), to the code of unit tests according to the configuration describing the syntax of the language, the testing framework and the interval library. 

The notation of ITL is easy to read and write and covers all notions of IEEE-1788 necessary to test an implementation. The following small example of test description from \cite{ITF1788Pub} gives a general understanding of ITL composition. 

\begin{verbatim}
/* Testing the addition function */
testcase addition.test {
    add [ -1.0, 1.0 ] [ empty ] = [ empty ];
    add [ 1.0, 2.0 ] [ 3.0, infinity ] = [ 4.0, infinity ];
    add [ 1.0, infinity ] [ -infinity, 4.0 ] = [ entire ];
    // using hexadecimal notation
    add [ 0X1.FFFFFFFFFFFFP+0 ] [ 0X1.999999999999AP-4 ] =
        [ 0X1.0CCCCCCCCCCC4P+1, 0X1.0CCCCCCCCCCC5P+1 ];
}
/* Testing the division function */
testcase division.test {
    div [ empty ] [ empty ] = [ empty ];
    div [ -30.0, 15.0 ] [ entire ] = [ entire ];
}
\end{verbatim}

ITF1788 inputs ITL-files and converts tests to the code of unit tests for specified language, test framework and library.

One can easily customize ITF1788 to support new programming languages, test frameworks and library specific features. The developer simply has to prepare several YAML configuration files. Additional flexibility can be obtained using implementation of Python callback functions for proper modifications of ITF1788 output.

Original ITF1788 engine and accompanying test sets for IEEE-1788 compliance testing have been developed by M. Kiesner, M. Nehmeier and J. Wolff von Gudenberg~\cite{ITF1788Pub}. Later O. Heimlich contributed a lot to the original project and to its own fork of ITF1788 \cite{ITF1788Codes}. The latter now contains:
\begin{itemize}
    \item configurations for programming languages C++, Octave, Julia, Python;
    \item configurations for test frameworks Boost Test Library, test frameworks for Octave, Julia and Python);
    \item plugins for libraries: 
        \begin{itemize}
            \item \texttt{libieeep1788, C-XSC, GAOL, Ibex} (in C/C++); 
            \item \texttt{interval package, INTLAB toolbox} (in GNU Octave); 
            \item \texttt{pyIbex} (in Python); 
            \item \texttt{ValidatedNumerics} (in Julia).
        \end{itemize}   
    \item almost 10,000 tests in ITL.
\end{itemize}   
Most of the tests accumulated in ITF1788 are derived from unit tests of C-XSC, FILIB, MPFI, libieeep1788.

\section{What is missing? A roadmap for testing}
\label{section-wishlist}

A first important requirement, which is perhaps too obvious but still
worth mentioning, is that a standard framework for tests must be open.
It should also be portable across architectures and languages, and be properly
documented.
The ITF1788 test suite has developed an interesting approach
to solve these particular issues. This is achieved by introducing a domain
specific language for storing the unit tests, and
allowing the user to pass the syntactic information necessary of the language
through the configuration of some files~\cite{ITF1788Pub}.

The unit tests of the ITF1788 suite cover numerous cases, including
some corner cases,
allow modularity and allow to test the tightness (equality assertion) or
correctness (asserting containment). Both of these aspects are desirable
and should be included with certain redundancy.
For instance, to avoid that an implementation of $\sin$ trivially returns
$[-1,1]$ and passes all containment assertions, tests that involve the
radius of the returned interval (to be less than two) can be introduced that
break the trivial implementation.

Yet, the unit tests of the ITF1788 suite are restricted to the {\tt binary64}
floating point format.
In our view, it is desirable to have the unit test coverage enlarged to include
unit tests for quadruple and extended precision, including
different precision values (e.g., 256, 512, 1024).
Similarly, concrete tests should be included for other IEEE754 numeric
formats, such as {\tt binary32}.

Other unit tests should choose randomly (floating-point) values $x \in \mathbf{x}$
and evaluate that the inclusion $f([x,x]) \subseteq f(\mathbf{x})$ holds true.
While these tests are certainly redundant, in particular if the implementation is
correct, they serve to check values not explicitly covered in previous tests,
in particular worst-case accuracy.
Such tests may be particularly interesting during the early stages of development,
considering specially the power function and the transcendental functions,
in particular for extended precision intervals.

In other communities, also a priori interested in reliable computations, "friendly competitions" are organized, to compare what each tool offers: methods, limitations, new aspects, etc. The code of the tools must be publicly available, and the benchmarks must be reproducible. See for instance \url{https://cps-vo.org/group/ARCH/FriendlyCompetition}.
Organizing such a friendly competition would be a sane motivation for developers, to update and increase the depth of unit tests.

\section{Conclusions}
\label{section-conclusion}

This paper focused on testing strategies for interval arithmetic libraries. During the years, different libraries have adopted different design choices and as such comparison and reproducibility has been a challenge. The IEEE 1788-2015 standardization is a big step to solve this issues, but broad adoption of the standard is a slow process. Moreover, testing is critical in software development and, due to its rigorous nature, even more central in interval arithmetic. It is thus fundamental to have an extensive unified testing framework. In this work, we sketched a possible structure for such framework. Hopefully, developers of interval libraries will join forces and share experience to develop an extensive open-source cross platform testing framework.
\vspace*{-3mm}
%
%
%
%

\end{document}